\documentclass[10pt]{amsart}

\usepackage{amssymb}
\usepackage{amsmath}
\allowdisplaybreaks
\usepackage{amsfonts}
\usepackage{fancybox}
\usepackage{empheq}
\usepackage{mathtools}
\usepackage{hyperref}

\usepackage[table]{xcolor}

\usepackage{color}
\usepackage{longtable,booktabs}
\usepackage{calrsfs} 
\usepackage{tikz-cd}
\usepackage[all,cmtip]{xy}

\usepackage{todonotes}
\usepackage{cleveref}

\definecolor{LightCyan}{rgb}{0.88,1,1}
\definecolor{LightCyan1}{rgb}{0.80,1,1}
\definecolor{Gray}{gray}{0.9}
\definecolor{Gray1}{gray}{0.95}

\usepackage{tikz}

\usepackage{changes}
\definechangesauthor[color=orange,name={Aristides Kontogeorgis}]{AK}

\newcommand{\Z}{\mathbb{Z}}

\newcommand{\Q}{\mathbb{Q}}

\newtheorem{theorem}{Theorem}

\theoremstyle{definition}

\DeclareRobustCommand{\mybox}[2][gray!20]{%
\begin{tcolorbox}[   
        breakable,
        left=0pt,
        right=0pt,
        top=0pt,
        bottom=0pt,
        colback=#1,
        colframe=#1,
        width=\dimexpr\textwidth\relax, 
        enlarge left by=0mm,
        boxsep=5pt,
        arc=0pt,outer arc=0pt,
        ]
        #2
\end{tcolorbox}
}

 \usepackage[theorems,breakable]{tcolorbox}
\tcbuselibrary{theorems}

 \newtcbtheorem[number within=section]{impdef}{Definition}%
 {colback=yellow!5,colframe=LightCyan!,fonttitle=\bfseries,coltitle=black}{th}

\title[Alexander modules]{A non-commutative differential module approach to Alexander modules}

\date{\today}

\author[A. Kontogeorgis]{Aristides Kontogeorgis}
\address{Department of Mathematics, National and Kapodistrian  University of Athens
Pane\-pist\-imioupolis, 15784 Athens, Greece}
\email{kontogar@math.uoa.gr}

\author[P. Paramantzoglou]{Panagiotis Paramantzoglou }
\address{Department of Mathematics, National and Kapodistrian University of Athens
Pane\-pist\-imioupolis, 15784 Athens, Greece}
\email{pan\_par@math.uoa.gr}

\date \today

\makeatletter
\newcommand{\aprod}{\mathop{\operator@font \hbox{\Large$\ast$}}}
\makeatother

\begin{document}

\begin{abstract}
The theory of R. Crowell on derived modules is approached within the theory of non-commutative differential modules. We also seek  analogies to the theory of cotangent complex from differentials in the commutative ring setting. Finally we give examples motivated from the theory of Galois coverings of curves.  
\end{abstract}

\maketitle


%
%
%
\section{Introduction}

In \cite{CrowellDerivedModule} R. Crowell defined the derived module corresponding to a homorphism $\psi:G \rightarrow H$ where $G$ is an arbitrary, in general non abelian group and $H$ is a second group, which in many interesting cases is assumed to be abelian. 
A $\psi$-derivation is a map such that for every 
$g_1,g_2\in G$ we have 
\begin{equation}
\label{psideriva}
\partial (g_1 g_2)=\partial(g_1)+\psi(g_1) \partial(g_2).
\end{equation}
The derived module $\mathcal{A}_\psi$, also known as the  $\psi$-differential module or as the Alexander module,  is the quotient module of the left $\Z[H]$-module $\bigoplus_{g\in G} \Z[H] dg$, generated by the symbols $dg$, for $g\in G$, 
divided by the left $\Z[H]$-module generated by elements of the form $d(g_1 g_2)-dg_1-\psi(g_1) dg_2$ for all $g_1,g_2 \in G$. The module $\mathcal{A}_\psi$ satisfies the following universal property:
\mybox[gray!10]{
For any left  $\Z[H]$-module $A$  and any 
$\psi$-derivation $\partial:G \rightarrow A$,  there exists a unique
$\Z[H]$-homomorphism  $\phi:\mathcal{A}_\psi \rightarrow A$, such that the following diagram, is commutative
\[
\xymatrix{
    G \ar[r]^d \ar[rd]_\partial &  \mathcal{A}_\psi  \ar[d]^{\phi} \\
  & A
}
\]
}
Moreover, suppose that the group $G$ admits a presentation in terms of generators and relations as 
\[
G=\langle x_1,\ldots,x_r | R_1=\cdots=R_s=1 \rangle,
\] 
and denote by $\pi$ the natural epimorphism from the free group generated by $x_1,\ldots,x_r$ to $G$.
The module $\mathcal{A}_\psi$ admits a free resolution over $\Z[H]$, 
\begin{equation}
\label{freeresolution}
\xymatrix{
\Z[H]^s \ar[r]^{Q_\psi} & \Z[H]^r 
\ar[r] & 
\mathcal{A}_\psi \ar[r] & 0
}
\end{equation}
where $Q_\psi$ is the matrix given by 
\begin{equation}
\label{presQ}
Q_\psi=
\left(
(\psi\circ \pi) \left(
\frac{\partial R_i}{
    \partial x_j
} \right)
\right), \text{ for } 1 \leq i \leq s, 1\leq j \leq r
\end{equation}
and $\frac{\partial R_i}{\partial x_j}$ is the Fox derivative of $R_i$ with respect to $x_j$, see \cite{FoxI}, \cite[chap. 3]{BirmanBraids},\cite[chap. 8]{Morishita2011-yw}.

The motivation of Crowell was in the theory of knots especially in the study of the link group in terms of the Wirtinger representation. As a matter of fact, the Alexander polynomial of a link can be expressed in terms of the fitting ideal corresponding to the free resolution given above, see \cite[chap. 9]{Morishita2011-yw}. The theory of $\psi$-differential modules can be extended to the case of pro-finite groups and in this form has interesting applications to Iwasawa theory using the ``Arithmetic topology'' view of point, see \cite[chap. 10,11]{Morishita2011-yw}.
In this setting it can also be used in order to define Galois representations, of the absolute Galois group $\mathrm{Gal}(\bar{\Q}/\Q)$, where the absolute Galois group is seen as a pro-finite analogon of the braid group, see \cite{Ihara1985-it}, \cite{MorishitaATIT}.


In commutative algebra there is a very well developed theory of differentials, see \cite[chap. 16]{Eisenbud:95},  which is an essential tool of modern Algebraic Geometry, see for example \cite[chap. II.8]{Hartshorne:77}. In this theory the module of differentials can be seen as an object representing the functor of derivations defined on the category of rings. Given a short exact sequence of rings, the relative cotangent sequence can be defined, see \cite[prop. 16.2]{Eisenbud:95}. Moreover, although the category of rings is not a abelian category a theory of higher cotangent functors can be developed, known as Andr\'{e}-Quillen homology, see \cite{MR0257068}, \cite{MR0352220}, \cite{MR2355775}.

In this article we would like to see the theory of $\psi$-differential modules in a similar setting. We attempt to use two variations of categories, namely the category $\mathcal{C}_1$ of short exact sequences of groups and the category $\mathcal{C}_2$ of short exact sequences, where the middle group is fixed. These two categories are explained in sections \ref{sec:catC1} and \ref{sec:catC2} respectively.

For a non-commutative ring $A$, derivations $d:A\rightarrow M$ are defined when $M$ is an $(A,A)$-bimodule, since the definition of a derivation $d:A\rightarrow M$ satisfies $d(xy)=xd(y)+d(x)y$. A similar construction as in the commutative case can be formed, see \cite[III.10 p. 567]{MR1727844}, \cite{MR1851201}.
\mybox[blue!7]{
\begin{theorem}
Consider  a continuous homomorphism of pro-$\ell$ groups $\psi:G\rightarrow H$, (a homomorphism $\psi:G\rightarrow H$ of finitely presented groups respectively) and define the ring $\mathbf{A}=\Z_\ell[[G]]$ ($\mathbf{A}=\Z[G]$ respectively). Define the  category of $(\mathbf{A},\mathbf{A})$-bimodules, where the action $x,y \in \mathbf{A}$ on $M$ is given by $x\cdot m \cdot y= \psi(x)\cdot m$. The non-commutative $(\mathbf{A},\mathbf{A})$-bimodule of differentials, which represents derivations coincides with Alexander module. 

Let $\mathcal{C}_1,\mathcal{C}_2$ be the categories of short exact sequences and short exact sequences with middle group fixed respectively. The Alexander module defines a functor from $\mathcal{C}_i$, $i=1,2$ the the category of $\Z_\ell$-modules ($\Z$-modules in the discrete group case), such that it sends epimorphisms to epimorphisms in the $\mathcal{C}_1$ case and is right exact in the $\mathcal{C}_2$-case.
\end{theorem}
}

Our interest in Alexander modules is motivated by the following geometric setting. 
Consider a Galois covering $\pi:\bar{Y}\rightarrow \mathbb{P}^1_{\bar{Q}}$ of the projective line ramified above a finite set of  points $S$,  $S\subset \mathbb{P}^1_{\bar{\Q}}$. When we extend the scalars from $\bar{\Q}$ to $\mathbb{C}$ we can  see $\bar{Y}$  as a compact Riemann surface of genus  $g\geq 2$. 
 The curve $Y_0=\bar{Y}-\pi^{-1}(S)$ is a topological covering of $X_s=\mathbb{P}^1_{\mathbb{C}}-S$, which 
can be described in terms of covering theory and corresponds to a subgroup $R_0$ of $\pi_1(X_s)$. 
 In general the group $R_0$ can be described using the Reidemeister-Schreier method, see \cite{ParamPartA19}. 
 Let $\Gamma$ be the closure of the  subgroup of $\mathfrak{F}_{s-1}$ generated by the stabilizers of ramification points, that is 
\begin{equation}
\label{GammaDef}
\Gamma=
  \langle x_1^{e_1},\ldots,x_{s}^{e_s} \rangle,
\end{equation}
where $e_1,\ldots,e_s$ are the ramification indices of the ramification points of $\pi:\bar{Y}\rightarrow \mathbb{P}^1$. 
The group $R=R_0/R_0 \cap \Gamma$ corresponds to the closed curve $\bar{Y}$ as a quotient of the hyperbolic plane. Our motivation was to understand the  homology group $H_1(\bar{Y},\Z)$ as a $\mathrm{Gal}(\bar{Y}/\mathbb{P}^1)$-module. This problem is essentially the dual problem of determining the Galois module structure of spaces of holomoprhic differentials as $\mathrm{Gal}(\bar{Y}/\mathbb{P}^1)$-modules, but here the study falls within the theory of integral representations, since $H_1(\bar{Y},\Z)$ is a free $\Z$-module, see \cite{ParamPartA19}. 

Also the braid group $B_{s-1}$, which is the mapping class group of the projective line with $s$ points removed is known to act on $H_1(\bar{Y},\Z)$ and in the same spirit (using the Arithmetic topology analogy) the absolute Galois group acts on $H_1(\bar{Y},\Z_\ell)$, see \cite{Ihara1985-it}, \cite{IharaCruz}, \cite{MorishitaATIT}. The later approach requires replacing the usual fundamental group by the pro-$\ell$ \'etale fundamental group, which in this case is the profinite completion of the usual fundamental groups. Having the second case in mind, we also consider the pro-$\ell$ versions of the Alexander module, see \cite[sec. 9.3]{Morishita2011-yw}.

This geometric situation can be expressed in terms of the short exact sequence of (pro-$\ell$) groups
\begin{equation}
\label{short-def}
1 \rightarrow  R=
\frac{R_0}{\Gamma \cap R_0} 
\cong 
\frac{R_0\cdot \Gamma}{\Gamma}
\rightarrow  
\frac{\mathfrak{F}_{s-1}}{\Gamma} 
\stackrel{\psi}{\longrightarrow} 
\frac{\mathfrak{F}_{s-1}}{R_0\cdot \Gamma}
 \rightarrow 1. 
\end{equation}
We define the ring 
\[
\mathcal{A}^{R_0,\Gamma}=\Z_\ell[[\frac{\mathfrak{F}_{s-1}}{R_0 \cdot \Gamma}]].
\]
If we assume that $\mathfrak{F}_{s-1}' \subset R_0$ then we the ring $\mathcal{A}^{R_0,\Gamma}$
is a commutative ring, and if all $e_i>1$,  then it is a group algebra corresponding to a finite commutative group $H=\mathfrak{F}_{s-1}/R_0 \cdot \Gamma$, that is 
$\mathcal{A}^{R_0,\Gamma}=\Z_\ell[H]$. 
In order to make the dependence clear we will denote the Alexander module in this setting by $\mathcal{A}_\psi^{R_0,\Gamma}$ instead of $\mathcal{A}_\psi$.
The Crowell exact sequence gives us that 
\cite[sec. 9.2, sec. 9.4]{Morishita2011-yw},
\begin{equation} \label{CrowellEx}
0 \rightarrow 
\left(
R
\right)^{\mathrm{ab}}
 =R/R'
\stackrel{\theta_1}{\longrightarrow} 
\mathcal{A}_\psi^{R_0,\Gamma}
\stackrel{\theta_2}{\longrightarrow}
\mathcal{A}^{R_0,\Gamma}
\stackrel{\varepsilon_{\mathcal{A}}}{\longrightarrow} 
\mathbb{Z}_\ell
\rightarrow 
0,
\end{equation}
where $R^{\mathrm{ab}}$ can be identified as the homology of the complete curve $\bar{Y}$.

In this way we have the group $R^{\mathrm{ab}}$ in a sequence of well understood $\Z[H]$-modules and this construction provides us with information on the $\Z[H]$-module structure of the first homology group.

In our study we have to see how  Alexander modules
corresponding to certain groups are related after taking quotients and for this an analogon of the   cotangent exact sequence is needed.

The Alexander module when $G$ is the free group $\mathfrak{F}_{s-1}$ and $H$ is the commutator group $\mathfrak{F}_{s-1}/\mathfrak{F}_{s-1}'$ is a free module on the ring $\mathcal{A}=\Z_\ell[[x_1,\ldots,x_{s-1}]]$ . In order to pass from the above case to the study of Alexander module for the map $\mathfrak{F}_{s-1}/\Gamma \rightarrow \mathfrak{F}_{s-1} / \Gamma \cdot \mathfrak{F}_{s-1}'$ requires relating objects in the category $\mathcal{C}_1$. This computation is explained in section \ref{sec:CatC1}. The transfer from $\mathfrak{F}_{s-1}$ to $\mathfrak{F}_{s-1}/\Gamma$ is important in the study of Braid group and $\mathrm{Gal}(\bar{\Q}/\Q)$-representations since it is related to the Gassner representation, see \cite[chap. 3]{BirmanBraids}, \cite{Ihara1985-it}.

On the other hand the passage from the Alaxander module for the map  $\mathfrak{F}_{s-1}/\Gamma \rightarrow \mathfrak{F}_{s-1} / \Gamma \cdot \mathfrak{F}_{s-1}'$ to the Alexander module for the map 
$\mathfrak{F}_{s-1}/\Gamma \rightarrow \mathfrak{F}_{s-1} / \Gamma \cdot R_0$ for some group $R_0$ satisfying $\mathfrak{F}_{s-1}' \subset R_0 \subset \mathfrak{F}_{s-1}$ requires working in the category $\mathcal{C}_2$, and this construction is explained in section \ref{relatingAlexMod}.

\subsection{Noncommutative modules of differentials}
We will fit the theory of derived modules  within the theory of derivations of noncommutative rings, see \cite[III.10 p. 567]{MR1727844}, \cite{MR1851201}. 

Consider the map of groups $\psi:G \rightarrow H$. The groups $G,H$ can be discrete finitely presented groups or pro-$\ell$ completions of discrete finitely presented groups. In the second case the map $\psi$ is assumed to be continuous. 
Define the noncommutative ring $\mathbf{A}=\Z_\ell[[G]]$ (resp. $\mathbf{A}=\Z[G]$ in the case of descrete groups). In order to make the presentation simpler we will write down the pro-$\ell$ case since the discrete case is exactly the same, one has simply to replace everywhere $\Z_\ell$ by $\Z$ and the completed group algebra $\Z_\ell[[G]]$ by the discrete group algebra $\Z[G]$.

 We begin by  defining the category of $(\mathbf{A},\mathbf{A})$-bimodules, where the action 
of $x,y\in \mathbf{A}$ on $M$ is given by 
\[
x\cdot m \cdot y = \psi(x) \cdot m, \text{ for all } x,y\in \mathbf{A}, m\in M.
\]
This means that $\mathbf{A}$ acts on $M$ from the left as a $\Z_\ell[[H]]$-module through the map $\psi$ and trivially on the right. So a derivation is a $\Z_\ell$-linear map with the additional property 
\[
d(xy)=\psi(x) d(y)+d(x)\cdot y= \psi(x) d(y)+ d(x), \qquad x,y\in \mathbf{A}.
\]
Of course this definition matches the definition of the $\psi$-derivation given in eq. (\ref{psideriva}). 
We can now consider the multiplication map
\[
m:\mathbf{A} \otimes_{\Z_\ell} \mathbf{A}
\rightarrow \mathbf{A},   m(x\otimes y)=xy.
\]
The module $\mathbf{A} \otimes_{\Z_\ell} \mathbf{A}$ is a $(\mathbf{A},\mathbf{A})$-bimodule with 
$\psi$-action on the left and trivial action on the right.

Following the exposition of Bourbaki we can consider the mapping 
\begin{equation}
\label{delta-map}
\delta: x\mapsto x(1\otimes 1)-(1\otimes1)x=x \otimes 1 -1\otimes 1,
\end{equation}
which is a derivation, i.e. $\delta(xy)=\delta(x)+\psi(x)\delta(y)$. The kernel of $\delta$ is known to be (see \cite[lemma 1]{MR1727844}) generated, as left $\mathbf{A}$-module, by the image of $\delta$, and also represents the space of derivations, in our bimodule setting, that is
\[
\mathrm{Hom}_{(\mathbf{A},\mathbf{A})}(
\mathrm{ker}\delta, M)\stackrel{\cong}{\longrightarrow} \mathrm{Der}_{\Z_\ell}(\mathbf{A},M), 
f\mapsto f\circ \delta. 
\]
By the form of eq. (\ref{delta-map}) we see that ($I_{\mathbf{A}}=I_{\Z_\ell[[G]]}$ is the augmentation module)
\[
\mathrm{ker}\delta = \mathbf{A}\delta(\mathbf{A}) =\Z_\ell[[H]]\otimes_\mathbf{A} I_{\mathbf{A}}, 
\]
which is in accordance to the computation in \cite[prop. 9.4]{Morishita2011-yw}.

Therefore the Alexander module corresponds to the object representing derivations in this setting and can be considered as a non-commutative space of differentials. In this point of view the Crowell sequence \cite[sec. 9.2, sec. 9.4]{Morishita2011-yw} is the analogon of the conormal  sequence in the commutative case, see \cite[prop. 16.3]{Eisenbud:95}.

Following this similarity we observe also 
that  the computation of the Alexander module given in 
eq. (\ref{presQ}) is the analogon of the conormal sequence given in section 16.1 in \cite[p. 387]{Eisenbud:95}.

\subsection{On the category $\mathcal{C}_1$ of short exact sequences of groups}

\label{sec:catC1}

The category of short exact sequences of groups, which has as objects short exact sequences and a homomorphism $f:A\rightarrow B$ of two short exact sequences $A,B$ is a map that makes the following diagram commutative:
\[
\xymatrix{
A \ar[d]^{f} &  1\ar[r] &  N_A \ar[r] \ar[d]^{f_N} & G_A \ar[r] 
  \ar[d]^{f_G} & H_A  \ar[d]^{f_H} \ar[r] & 1 
\\
B  &  1 \ar[r] &  N_B \ar[r] & G_B \ar[r] & H_B \ar[r] & 1
}
\]
The Alexander module  functor sends the sequence $A$ to the Alexander module $\mathcal{A}_{G_A\rightarrow H_A}$ which will be denoted by 
$\mathcal{A}_A$. The map $f=(f_N,f_G,f_H)$ from $A\rightarrow B$ induces a map 
\[
\mathcal{A}(f)=f_H \otimes f_G: 
\Z_\ell[[H_A]] 
\otimes_{\Z_\ell[G_A]}
I_{\Z_\ell[G_A]}
\longrightarrow
\Z_\ell[[H_B]] 
\otimes_{\Z_\ell[G_B]}
I_{\Z_\ell[G_B]}.
\]
Several problems in this attempt arrise: The category of short exact sequences is not an abelian category, since the category of groups is not 
an abelian category (it is not even additive). 
For example the category of groups does not have cokernels in the category, since the image of a group homomorphism is not necessary normal. However, as one does in non-abelian cohomology theory  \cite[Appendix]{SeL},
in the category of pointed sets cokernels exist and then the snake lemma can provide us with a cokernel in the category of pointed sets, see the discussion in 
\cite{MathOver53249}.

Assume now that we have a short exact sequence of short exact sequences
\begin{equation}
\label{sesofses}
0\rightarrow A \rightarrow B \rightarrow C \rightarrow 0.
\end{equation}
We can prove that $\mathcal{A}_B\rightarrow \mathcal{A}_C$ is onto. Indeed, by propositions 2 and 6  in \cite[II.3,p.245-252]{MR1727844} the composite map
$\mathcal{A}_{B\rightarrow C}$
\[
\xymatrix@C-1pc{
\mathcal{A}_B \ar@{=}[d] 
\ar[rr]^{\mathcal{A}_{B\rightarrow C}}
& & \mathcal{A}_C \ar@{=}[d]\\
\Z_\ell[[H_B]]
\otimes_{\Z_\ell[[G_B]]}
I_{\Z_\ell[G_B]}
\ar[r] 
& 
\Z_\ell[[H_C]]
\otimes_{\Z_\ell[[G_B]]}
I_{\Z_\ell[G_C]}
\ar[r]
&
\Z_\ell[[H_C]]
\otimes_{\Z_\ell[[G_C]]}
I_{\Z_\ell[G_C]}
}
\]
is onto. Again by the description of kernel given in proposition 6 it seems that we don't have exactness at middle sequence $B$ in the short exact sequence of short exact sequences given in eq. (\ref{sesofses}).

\subsubsection{An example:Alexander module for the commutator group}

\label{sec:CatC1}
%
%


Let us now consider a special case of the above construction where $G_B$ is the the pro-$\ell$ free group in $s-1$ generators $\mathfrak{F}_{s-1}$, $N_B=\mathfrak{F}_{s-1}'$ and $H_B=\mathfrak{F}_{s-1}/\mathfrak{F}_{s-1}'$, while $G_C=\mathfrak{F}_{s-1}/\Gamma$, $N_C=\mathfrak{F}_{s-1}'\cdot \Gamma/\Gamma$, $H_C=\mathfrak{F}_{s-1}/\mathfrak{F}'_{s-1}\cdot \Gamma$ and $\Gamma$ is a normal closed subgroup of $\mathfrak{F}_{s-1}$ generated by $r$-elements. So we have the following map of short exact sequences:
\[
\xymatrix{
B \ar[d]^{f} &  
1\ar[r] & \mathfrak{F}_{s-1}'  \ar[r] \ar[d] &
 \mathfrak{F}_{s-1} \ar[r] 
  \ar[d]
  & \frac{\mathfrak{F}_{s-1}}{\mathfrak{F}_{s-1}'}  \ar[d] \ar[r] & 1 
\\
C  &  1 \ar[r] & \frac{\mathfrak{F}_{s-1}' \cdot \Gamma}{\Gamma} \ar[r]  &
\frac{\mathfrak{F}_{s-1}}{\Gamma} \ar[r]  & 
\frac{\mathfrak{F}_{s-1}}{\mathfrak{F}'_{s-1} \cdot \Gamma}
 \ar[r]  &
1
}
\]

 Motivated by the applications of the differential module to the computations of homology groups of coverings of the projective line, we will use this construction for the group $\Gamma$ given in eq. (\ref{GammaDef}) but this construction can also  be used  for a more general group $\Gamma$.
The Alexander module corresponding to the last line of the above diagram is  a module over the ring 
\[
\mathcal{A}^{\mathfrak{F}'_{s-1},\Gamma}=\Z_\ell[[\mathfrak{F}_{s-1}/\mathfrak{F}'_{s-1} \cdot \Gamma]].
\] 
while the Alexander module corresponding to the top line is a module over the ring 
\[
\mathcal{A}^{\mathfrak{F}_{s-1}',\{1\}}=\Z_\ell[[\mathfrak{F}_{s-1}^{\mathrm{ab}}]]\cong\Z_\ell[[u_1,\ldots,u_{s-1}]]=:\mathcal{A}.
\]
We also consider the  free resolution of Alexander modules as given in eq. (\ref{freeresolution}). 
In the first row we consider the group $\mathfrak{F}_{s-1}$ as the quotient of the free pro-$\ell$ group is $s$-generators modulo the relation $x_1x_2\cdots x_s=1$. In the second row the group $\mathfrak{F}_{s-1}/\Gamma$ is considered as the quotient of the free group in $s$-generators modulo the relation $x_1x_2\cdots x_s=1$ and the $r$-relations generating $\Gamma$. 
\begin{equation}
\label{presentFree}
\xymatrix@C=15pt{
\mathcal{A} \ar[r]^-{Q_1}
&
\mathcal{A}^s \ar[r]^-{\psi_1}  \ar[d]^{\phi_2}&
\mathcal{A}^{s-1} \ar[r] \ar[d]^{\phi_3} & 0
\\
\left(\mathcal{A}^{\mathfrak{F}_{s-1}',\Gamma}\right)^{r+1} \ar[r]^{Q_2} &
\left(\mathcal{A}^{\mathfrak{F}_{s-1}',\Gamma}\right)^{s} \ar[r]^-{\psi_2} &
\mathcal{A}^{\mathfrak{F}'_{s-1},\Gamma}_\psi \ar[r] & 0
}
\end{equation}
where $Q_1,Q_2$ are the maps appearing in eq. (\ref{freeresolution}). In particular the map $Q_1$ sends 
\[
\mathcal{A} \ni \beta \mapsto \beta\cdot( 1,x_1,x_1x_2,\ldots,x_1\cdot x_2 \cdots x_{s-1} ).
\]
The vertical map $\phi_2$  is the  reduction modulo $\Gamma$ and it is onto.
The image $\phi_3(a)$ for $a\in \mathcal{A}^{s-1}$
 is defined by selecting $b\in \mathcal{A}^{s}$ such that $\psi_1(b)=a$, and 
then $\phi_3(a)=\psi_2\circ\phi_2(b)$ as seen in the diagram bellow:
\[
\xymatrix{
  b \ar[r]^{\psi_1} \ar[d]^{\phi_2} & a \ar[d]^{\phi_3} \\
  \phi_2(b) \ar[r]^-{\psi_2} & \phi_3(a)=\psi_2 \circ \phi_2(b)
}
\]
This definition is independent from the selection of $b$.  
We have 
\begin{equation}
\label{kerf3}
\mathrm{ker}(\phi_3)=\psi_1 \big( \phi_2^{-1}  (\mathrm{Im}(Q_2))\big).
\end{equation}
For the commutator group of a quotient ($B\lhd A$) we have
$
(A/B)'=A' B/B
$ so 
\[
\left(\frac{
\mathfrak{F}_{s-1}'\cdot \Gamma
}
{
\Gamma
}
\right)' 
\!\!=\!\!
\left(\frac{
\mathfrak{F}_{s-1}'
}
{
\Gamma \cap \mathfrak{F}_{s-1}'
}
\right)'
\!\!=\!\!
\frac{
\mathfrak{F}_{s-1}'' (\Gamma \cap \mathfrak{F}_{s-1}')
}
{
\Gamma \cap \mathfrak{F}_{s-1}'
}
\!\!=\!\!
\frac{
\mathfrak{F}_{s-1}'' 
}
{
\Gamma \cap \mathfrak{F}_{s-1}' \cap \mathfrak{F}_{s-1}''
}
\!\!=\!\!
\frac{
\mathfrak{F}_{s-1}'' 
}
{
\Gamma \cap \mathfrak{F}_{s-1}''
}
\!\!=\!\!
\frac{
\mathfrak{F}_{s-1}'' \Gamma
}
{
\Gamma
}.
\]
We finally have
\begin{equation}
\label{*p9}
\left(  \left( 
\frac{
  \mathfrak{F}_{s-1}
}{
  \Gamma
} \right)' \right)^{\mathrm{ab}}
=
\left(\frac{
\mathfrak{F}_{s-1}'\cdot \Gamma
}
{
\Gamma
}
\right)^{\mathrm{ab}} 
=
\frac{\mathfrak{F}_{s-1}'\cdot \Gamma}{\mathfrak{F}_{s-1}''\cdot \Gamma}.
\end{equation}
The corresponding Crowell sequences are given:
\begin{equation}
\label{crowGamma}
\xymatrix{
0 \ar[r] & 
\frac{
\mathfrak{F}_{s-1}'
}
{
\mathfrak{F}_{s-1}''
} 
\ar[r] \ar[d]^\phi &
 \mathcal{A}^{s-1} \ar[r] \ar[d]^{\phi_3} &
 \mathcal{A} \ar[r] \ar[d]^{\mathrm{mod}\Gamma} & \mathbb{Z}_\ell
\ar[r] \ar@{=}[d] & 0 
 \\
 0 \ar[r] & 
\left(
\frac{
\mathfrak{F}_{s-1}'\cdot \Gamma
}
{
\Gamma
}
\right)^{\mathrm{ab}} 
=
\frac{\mathfrak{F}_{s-1}'\cdot \Gamma}{\mathfrak{F}_{s-1}''\cdot \Gamma}
\ar[r]  &
 \mathcal{A}^{\mathfrak{F}_{s-1}',\Gamma}_\psi \ar[r]  &
 \mathcal{A}^{\mathfrak{F}_{s-1}',\Gamma} \ar[r] &\Z_\ell \ar[r]  & 0
}
\end{equation}
We have 
\[
\mathrm{ker}\phi=
\frac{\mathfrak{F}_{s-1}''\cdot\Gamma \cap \mathfrak{F}_{s-1}'}
{\mathfrak{F}_{s-1}''}.
\]
Also by the sharp five lemma, or the diagram in eq. (\ref{presentFree}) the map $\phi_3$ is onto, see \cite{nlab:five_lemma}. 
\subsection{On the category $\mathcal{C}_2$ of short exact sequences with middle group fixed}
\label{sec:catC2}
Consider now the category of shortexact sequences with middle group $G$ fixed and maps $(f_N,\mathrm{id},f_H)$. Essentially this is the category with objects the pairs $A=(H_A,f_A)$, $f_A:G\rightarrow H_A$
and functions  $f_{A\rightarrow B}:H_A \rightarrow H_B$, so that $  f_{A\rightarrow B} \circ f_A =f_B$, that is the following diagram is commutative:
\[
\xymatrix{
  G \ar[rd]_{f_B} \ar[r]^{f_A} & H_A \ar[d]^{f_{A\rightarrow B}} \\
  &  H_B
}
\]
For this category we can define again the 
Alexander module functor sending the pair $(H,f)$ to the Alexander module $\mathcal{A}_f$. 

In particular it is very interesting to assume that all groups $H_A$ are commutative. In this case cokernels exist. Moreover for a short exact sequence of pairs
\[
0 \rightarrow A \rightarrow B \rightarrow C \rightarrow 0
\]
we have that the corresponding   sequence of Alexander modules is exact:
\[
\xymatrix@C-1.5pc{
\mathcal{A}_A \ar@{=}[d] \ar[r] 
&
\mathcal{A}_B  \ar@{=}[d] \ar[r]
&
\mathcal{A}_C \ar@{=}[d] \ar[r] 
& 0 
\\
\Z_{\ell}[[H_A]]\otimes_{\Z_\ell[[G]]} I_{\Z_\ell[[G]]}
\ar[r] &
\Z_{\ell}[[H_B]]\otimes_{\Z_\ell[[G]]} I_{\Z_\ell[[G]]}
\ar[r] &
\Z_{\ell}[[H_C]]\otimes_{\Z_\ell[[G]]} I_{\Z_\ell[[G]]}
\ar[r] &
0
\\ 
}
\]
by using proposition 5 in  \cite[II.3,p.245-252]{MR1727844}.

\subsubsection{Relating Alexander modules in the category $\mathcal{C}_2$}
\label{relatingAlexMod}
Assume that $\mathfrak{F}_{s-1}' \subset R_0 \subset \mathfrak{F}_{s-1}$. To the commutative diagram
\[
\xymatrix{
  1 \ar[r] & R=\frac{R_0}{R_0 \cap \Gamma} \ar[r] & \mathfrak{F}_{s-1}/\Gamma \ar[r]^-{\psi} & 
  \mathfrak{F}_{s-1}/R_0 \cdot \Gamma \ar[r] & 1 
  \\
  1 \ar[r] & \;\;\;
  \frac{\mathfrak{F}_{s-1}' }
  {
  \mathfrak{F}_{s-1}' \cap \Gamma
  }
  \ar[r] \ar@{^{(}->}[u] & \mathfrak{F}_{s-1}/\Gamma \ar[r]^-{\psi} \ar@{=}[u]& 
  \mathfrak{F}_{s-1}/\mathfrak{F}_{s-1}' \cdot \Gamma \ar[r] \ar[u] & 1
}
\]
we can attach two related  Crowell sequences:
\[
\xymatrix{
  0 \ar[r] & \ar[r] R^{\mathrm{ab}} & \mathcal{A}^{R_0,\Gamma}_\psi \ar[r]& \mathcal{A}^{R_0,\Gamma} \ar[r] & \mathbb{Z}_\ell \ar[r] &  0
  \\
  0 \ar[r] & \ar[r] \ar@{^{(}->}[u]^{\theta_1} 
\frac{
  \mathfrak{F}'_{s-1}\cdot \Gamma
  }
  {
  \mathfrak{F}_{s-1}'' \cdot \Gamma
  }
  & \mathcal{A}^{\mathfrak{F}_{s-1}',\Gamma}_\psi \ar[r] \ar[u]^{\theta_2} 
  & \mathcal{A}^{\mathfrak{F}'_{s-1},\Gamma} \ar[r] \ar[u]
  & \mathbb{Z}_\ell \ar[r] \ar@{=}[u] &  0
}
\]
Using the same argument as in eq. (\ref{*p9}) we have that $R^{\mathrm{ab}}=R_0'\cdot \Gamma / R_0'' \cdot \Gamma$.
The map $\theta_1$ is well defined with kernel $R_0'' \cdot \Gamma/\mathfrak{F}_{s-1}''\cdot \Gamma$. The map $\theta_2$ is defined from the two corresponding free resolutions:
\[
\xymatrix{
 \Z_\ell[[\mathfrak{F}_{s-1}/R\cdot\Gamma]]^{s+1} \ar[r]^{Q_1} & 
 \Z_\ell[[\mathfrak{F}_{s-1}/R\cdot\Gamma]]^s   \ar[r]^-{\pi_1} &
 \mathcal{A}^{R,\Gamma}_\psi \ar[r] & 0 
 \\
 \Z_\ell[[\mathfrak{F}_{s-1}/\mathfrak{F}_{s-1}'\cdot\Gamma]]^{s+1} \ar[r]^{Q_2} \ar[u]^{\phi_1}
 & \Z_\ell[[\mathfrak{F}_{s-1}/\mathfrak{F}_{s-1}'\cdot\Gamma]]^s   \ar[r]^-{\pi_2} \ar[u]^{\phi_2} &
 \mathcal{A}^{\mathfrak{F}_{s-1}',\Gamma}_\psi  \ar[u]^{\theta_2} \ar[r] & 0  
}
\] 
Indeed, for $a \in \mathcal{A}^{\mathfrak{F}_{s-1}',\Gamma}$ we select any
$b \in \Z_\ell[[\mathfrak{F}_{s-1}/\mathfrak{F}_{s-1}'\cdot\Gamma]]^s $ and we set
\[
\theta_2(a)=\pi_1\circ\phi_2(b).
\]
Then $\theta_2$ is well defined, i.e. independent from the selection of $b$.  The kernel of $\theta_2$ is
\[
 \mathrm{ker}(\theta_2)=\pi_2
 \left(
 \phi_2^{-1}(\mathrm{Im}(Q_1))
 \right).
\]

Observe also that both $\mathcal{A}^{R_0,\Gamma}_\psi$ and $\mathcal{A}^{\mathfrak{F}'_{s-1},\Gamma}_\psi$ are the cokernel of the same set of equations since the matrix $Q$ depends only on the quotient $\mathfrak{F}_{s-1}/\Gamma$. The difference is that they are modules over different rings. 

 \def\cprime{$'$}

\end{document}